\documentclass[12pt]{amsart}
\usepackage[utf8]{inputenc}
\usepackage{color}
\usepackage{times}
\usepackage[hidelinks]{hyperref}
\usepackage{enumerate,latexsym}
\usepackage{latexsym}
\usepackage{subcaption}

\usepackage{amsmath,amsthm,amsfonts,amssymb}
\usepackage{graphicx}
\usepackage{dsfont}

\usepackage{tikz}
\usetikzlibrary{arrows.meta}
\usetikzlibrary{knots}
\usetikzlibrary{hobby}
\usetikzlibrary{arrows,decorations,decorations.markings}

\usetikzlibrary{fadings}

\makeatletter
\def\namedlabel#1#2{\begingroup
 #2%
 \def\@currentlabel{#2}%
 \phantomsection\label{#1}\endgroup
}
\makeatother

\usepackage[lite]{amsrefs}

\renewcommand{\PrintDOI}[1]{\href{http://dx.doi.org/\detokenize{#1}}{doi: \detokenize{#1}}%
	\IfEmptyBibField{pages}{, (to appear in print)}{}}

\theoremstyle{plain}
\newtheorem*{theorem*}{Theorem}
\newtheorem*{thmex*}{Theorem~\ref{example}}
\newtheorem*{thmasymp*}{Theorem~\ref{thmAsymp}}
\newtheorem{theorem}{Theorem}[section]

\newtheorem{proposition}[theorem]{Proposition}

\newtheorem{example}[theorem]{Example}

\theoremstyle{definition}
\newtheorem{definition}[theorem]{Definition}

\newcommand{\R}{\mathbb{R}}

\newcommand{\Z}{\mathbb{Z}}

\newcommand{\ben}{\begin{enumerate}}
\newcommand{\een}{\end{enumerate}}

\newcommand{\ed}{\end{document}}

\definecolor{rrr}{rgb}{.9,0,.1}

\definecolor{rr}{rgb}{.8,0,.3}

\usepackage{pdfsync}
\graphicspath{ {./images/} }

\title[Enhancement of the Coloring Invariant]{Enhancement of the Coloring Invariant for Folded Molecular Chains}

\author[J. Ceniceros]{Jose Ceniceros}
\address{Hamilton College, Clinton, NY, USA}
\email{jcenicer@hamilton.edu}

\author[M. elhamadi]{Mohamed Elhamdadi}
\address{University of South Florida, Tampa, Florida, USA}
\email{emohamed@usf.edu}

\author[A. Mashaghi]{Alireza Mashaghi}
\address{Leiden University, Leiden, The Netherlands}
\email{a.mashaghi.tabari@lacdr.leidenuniv.nl}

\begin{document}

\maketitle

\begin{abstract}

Folded linear molecular chains are ubiquitous in biology. Folding is mediated by intra-chain interactions that "glue" two or more regions of a chain. The resulting fold topology is widely believed to be a determinant of biomolecular properties and function. Recently, knot theory has been extended to describe the topology of folded linear chains such as proteins and nucleic acids. To classify and distinguish chain topologies, algebraic structure of quandles has been adapted and applied. However, the approach is limited as apparently distinct topologies may end up having the same number of colorings. Here, we enhance the resolving power of the quandle coloring approach by introducing Boltzmann weights. We demonstrate that the enhanced coloring invariants can distinguish fold topologies with an improved resolution. 

\end{abstract}
\tableofcontents

\section{Introduction}
Biomolecular chains such as proteins and nucleic acids fold into a large repertoire of conformations, with distinct roles in health and disease \cites{Dobson2002, Corces2018}. Resolving the topological complexity of the proteome provides insights into major fundamental and applied questions in biology, molecular medicine and bioinspired engineering \cite{Trends2020Scalvini}. It reveals the design principles of biomolecular systems and thus provides structural understandings of misfolded molecules. Furthermore, topological information may shed light on how folding is regulated in cells; molecular chaperones are known to guide folding of biomolecular chains, yet their working principles are not fully understood. It is speculated that they promote or suppress the formation of certain topologies \cite{ACS2020Heidari}. Molecular topology is also relevant to engineering fields. Inspired by the topological information, one can design folded molecular chains that mimic biomolecules or carry totally new functions with applications in medicine and engineering \cite{Trends2020Scalvini}. 

For nearly two centuries, knot theory has been developed and used to describe topology of folded chains. Knot theory has been widely applied to various areas of engineering and science; it has also been used to describe proteins and nucleic acids \cites{DS17, 7, mashaghi2014circuit}. However, most of structurally identified biomolecules such as proteins fall into the category of “unknot”, limiting the resolution of conventional knot theoretic approaches \cite{Sulkowska2012}. This is due to the fact that intra-chain interactions or bonds are key to folding of biomolecular chains, and knot theoretic approaches that exclusively focus on chain crossings will not be able to resolve topological features in such folded chains \cite{Trends2020Scalvini}. Recent extensions of knot theory addressed this shortcoming by including singularities in knot projections and expanding the Reidemeister move set \cite{ADEM}. The algebraic structure of quandles \cite{EN} can then be applied to distinguish the topologies of main protein chains. The approach provided a proof of concept that knot theory can be extended to provide a resolution beyond what was previously anticipated \cite{ADEM}. However, the current theory is still significantly limited in resolution. Chains whose knot projection have different “quandle coloring” are topologically distinct, however, identical quandle coloring can be seen in folded chains with distinct topologies. Here, we further extend the theoretical framework to increase the resolving power of the topological approach. We propose a protocol based on quandle coloring and Boltzmann weights for resolving topological complexity of folded linear chains.

The article is organized as follows.  In Section~\ref{RMKT}, we review the basics of knot theory.  Section~\ref{QSB} gives the diagrammatic needed to motivate the definitions of quandle, singquandle and bondle.  It also contains some examples.  In Section~\ref{EBCI}, we introduce an \emph{enhancement} of the bondle counting invariant using Boltzmann weights at crossings and bonds.  In Section~\ref{Examples}, we show the strength of the enhanced invariant by giving examples of folded molecular chains having the same bondle counting invariant but are distinguish by the enhanced invariant.   

\section{Knots, Singular Knots, and Folded Molecular Chains}\label{RMKT}

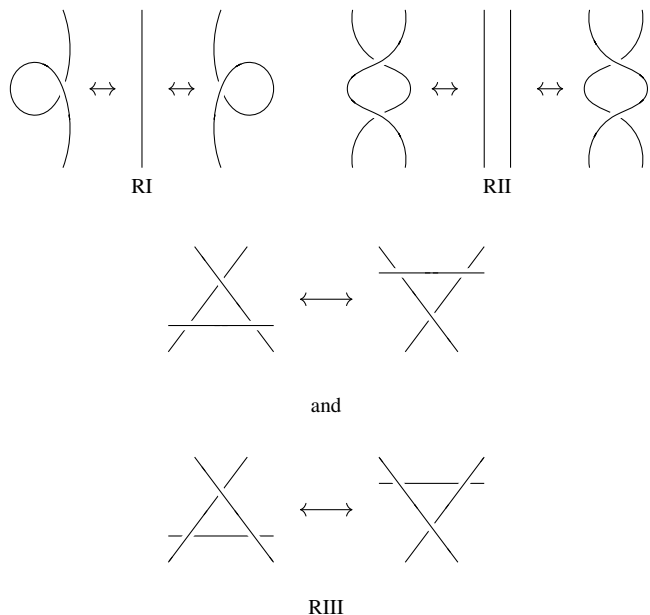
\begin{figure}[h]
\begin{tikzpicture}[scale=.7,use Hobby shortcut]
\begin{knot}[
consider self intersections=true,
ignore endpoint intersections=false,
]
\strand(-7,-1.5)..(-7,0)..(-7.5,.5)..(-8,0)..(-7.5,-.5)..(-7,0)..(-7,1.5); 
\draw[<->] (-6.5,0)..(-6,0);
\strand(-5.5,-1.5)..(-5.5,1.5);
\draw[<->] (-5,0)..(-4.5,0);
\strand(-4,-1.5)..(-4,0)..(-3.5,.5)..(-3,0)..(-3.5,-.5)..(-4,0)..(-4,1.5); 

\strand(-.5,-1.5)..(-1,-.5)..(-1.6,0)..(-1,.5)..(-.5,1.5);   
\strand(-1.5,-1.5)..(-1,-.5)..(-.4,0)..(-1,.5)..(-1.5,1.5);
\draw[<->] (0,0).. (.5,0);
\strand(1,-1.5)..(1,1.5);    
\strand(1.5,-1.5)..(1.5,1.5);
\draw[<->] (2,0)..(2.5,0);
\strand(3,-1.5)..(3.5,-.5)..(4.1,0)..(3.5,.5)..(3,1.5);    
\strand(4,-1.5)..(3.5,-.5)..(2.9,0)..(3.5,.5)..(4,1.5);

\strand(-4.5,-7)..(-3,-9);
\strand(-5,-9)..(-3.5,-7);  
\strand(-3,-8.5)..(-5,-8.5);
\draw[<->] (-2.5,-8)..(-1.5,-8);
\strand(-1,-7)..(.5,-9);
\strand(-.5,-9)..(1,-7); 
\strand(-1,-7.5)..(1,-7.5);

\strand(-3,-4.5)..(-5,-4.5);
\strand(-4.5,-3)..(-3,-5);
\strand(-5,-5)..(-3.5,-3);    
\draw[<->] (-2.5,-4)..(-1.5,-4);
\strand(1,-3.5)..(-1,-3.5);
\strand(-1,-3)..(.5,-5);
\strand(-.5,-5)..(1,-3); 
\end{knot}

\node[below] at (-5.5,-1.5) {\tiny RI};
\node[below] at (1.25,-1.5) {\tiny RII};
\node at (-2,-6) {\tiny and};
\node[below] at (-2,-9.5) {\tiny RIII};
\end{tikzpicture}
	\caption{Reidemeister Moves RI, RII, and RIII.}
		\label{reidemeister}
\end{figure}

A \emph{knot} is a closed curve in the $3$-dimensional space that does not intersect itself.  In other words, it is an embedding of $S^1,$ the circle, into three dimensions, called a \emph{conformation}, given by a function $g: [0,2\pi] \rightarrow \mathbb{R}^3$ where $g(0)= g(2\pi)$. We think of the knot as the image of this map $g$.  One of the central questions in knot theory is given two knot conformations, can one knot be deformed without passing through itself to obtain the other knot? If so, we say that the two knot conformations are equivalent. Therefore, a knot is an equivalence class of equivalent knot conformations. Although knots are 3-dimensional objects, it is convenient to focus on a knot conformation projection in a specific direction onto a plane, called a \emph{knot diagram}. Therefore, the question of knot equivalence can be reformulated in terms of knot diagrams. Given two knot diagrams, can one diagram be deformed to give the other diagram? If so, the two diagrams are said to represent the same knot.  In the 1920s, Kurt Reidemeister proved a theorem stating that two diagrams represent the same knot if one diagram can be obtained from the other by continuous deformation in the plane and the so-called \emph{Reidemeister} moves RI, RII and RIII as can be seen in Figure~\ref{reidemeister}.

\begin{figure}[h]
\begin{tikzpicture}[use Hobby shortcut,scale=1,add arrow/.style={postaction={decorate}, decoration={
  markings,
  mark=at position 1 with {\arrow[scale=1,>=stealth]{>}}}}]
\begin{knot}[
consider self intersections=true,
 ignore endpoint intersections=false,
]
    

\strand(0,-1)..(-.5,-.5)..(-1,0)..(-.5,.5)..(0,1);   
\strand(-1,-1)..(-.5,-.5)..(0,0)..(-.5,.5)..(-1,1);
\draw[<->] (.5,0).. (1,0);

\strand(1.5,-1)..(2,-.5)..(2.5,0)..(2,.5)..(1.5,1);    
\strand(2.5,-1)..(2,-.5)..(1.5,0)..(2,.5)..(2.5,1);


\strand(-.5,-3.5)..(-2.5,-3.5);
\strand(-2.5,-4)..(-1,-2);   
\strand(-2,-2)..(-.5,-4);
 \draw[<->] (.5,-3)..(1.5,-3);
\strand(4,-2.5)..(2,-2.5);
\strand(2.5,-4)..(4,-2); 
\strand(2,-2)..(3.5,-4);


\strand(-2.5,-7)..(-1,-5);   
\strand(-2,-5)..(-.5,-7);
\strand(-.5,-6.5)..(-2.5,-6.5);
 \draw[<->] (.5,-6)..(1.5,-6);
\strand(2.5,-7)..(4,-5); 
\strand(2,-5)..(3.5,-7);
\strand(4,-5.5)..(2,-5.5);
\end{knot}

\node[below] at (.75,0) {\tiny $sRII$};
\node[below] at (1,-3) {\tiny $sRIIIa$};
\node[below] at (1,-6) {\tiny $sRIIIb$};

\node[circle,draw=black, fill=black, inner sep=0pt,minimum size=5pt] (a) at (-.5,-.5) {};
\node[circle,draw=black, fill=black, inner sep=0pt,minimum size=5pt] (a) at (2,.5) {};

\node[circle,draw=black, fill=black, inner sep=0pt,minimum size=5pt] (a) at (-1.5,-2.7) {};
\node[circle,draw=black, fill=black, inner sep=0pt,minimum size=5pt] (a) at (3,-3.4) {};

\node[circle,draw=black, fill=black, inner sep=0pt,minimum size=5pt] (a) at (-1.5,-5.7) {};
\node[circle,draw=black, fill=black, inner sep=0pt,minimum size=5pt] (a) at (3,-6.4) {};
\end{tikzpicture}
	\caption{The Reidemeister moves involving singular crossings.}
		\label{SRmoves}
\end{figure}
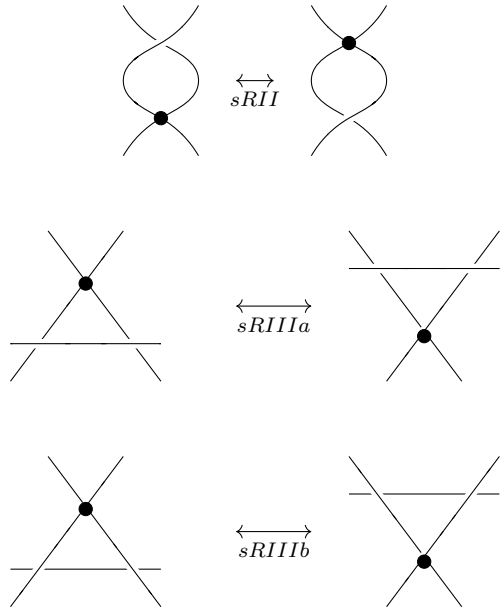

In knot theory, it is well known that any knot can be transformed to become the \emph{unknot} by a finite number of crossing changes \cite{Adams}.  If one thinks of a crossing change happening as time is changing, then at some point in time the two branches meet, thus making a point of self-intersection.  Extending knot theory by adding self-intersection points gives what is known as \emph{singular knot theory} \cite{Dancso}. Therefore, a singular knot is a knot, but with a finite number of transversal self-intersections. We call these self-intersections, \emph{singular crossings}. Similar to classical knots it is convenient to focus on diagrams of singular knots. Therefore, we will be projecting a singular knot in a particular direction onto a plane, to obtain a \emph{singular knot diagram}. In a singular knot diagram, we will define \emph{classical crossings} as a place were one strand goes under another and a \emph{singular crossing} will be defined as a place where two strands transversely intersect. Similar to classical knot theory a central question in singular knot theory, is given two singular links are the two topologically equivalent. Although the Reidemeister moves in Figure ~\ref{reidemeister} are able to completely describe the equivalence of classical knots, the addition of singular crossings requires three more Reidemeister moves, see Figure~\ref{SRmoves}. We will refer to the set consisting of the Reidemeister moves in Figure~\ref{reidemeister} and the singular Reidemeister moves in Figure~\ref{SRmoves}, as the \emph{generalized Reidemeister moves}. It is well known that two singular links $K_1$ and $K_2$ are equivalent if and only if there exists a sequence of generalized Reidemeister moves and continuous deformation in the plane that transforms a diagram of $K_1$ to a diagram of $K_2$.

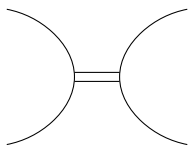
\begin{figure}[h]
\begin{tikzpicture}[use Hobby shortcut,scale=.6,add arrow/.style={postaction={decorate}, decoration={
  markings,
  mark=at position 1 with {\arrow[scale=1,>=stealth]{>}}}}]
\begin{knot}

 \strand (-3.5,1.5)..(-2.5,1)..(-2,0)..(-2.5,-1)..(-3.5,-1.5);
 \strand (.5,1.5)..(-.5,1)..(-1,0)..(-.5,-1) ..(.5,-1.5);

\end{knot}

 \draw (-2,.1) to (-1,.1);
 \draw (-2,-.1) to (-1,-.1);

\end{tikzpicture}
\vspace{.2in}
		\caption{Representation of a singular crossing in a molecular chain diagram.}
		\label{int}
\end{figure}

Following the conventions introduced in \cite{ADEM}, we will use singular knot theory to model folded molecular chains. We can think of a folded molecular chain as a conformation, given by a function $g:[0,1] \rightarrow \R^3$, which we will call a \emph{folded molecular chain model}. A folded molecular chain model differs from a singular knot since, in this case, we are considering a line segment with self-intersections in 3-space and not a closed curve with self-intersection in 3-space. It will be convenient for computational purposes to consider projections of a folded molecular chain model to obtain a \emph{folded molecular chain diagram}. We will define singular crossings to exist where a molecular chain has intra-chain interactions. Therefore, we will illustrate singular crossing in a folded molecular chain diagram by a small rectangle with two parallel edges on two strands, see Figure~\ref{int}. For a more detailed description of the model, we refer the reader to \cite{ADEM}.

\begin{figure}[h]
\begin{tikzpicture}[use Hobby shortcut,scale=.4]
\begin{knot}[
consider self intersections=true,
  ignore endpoint intersections=false,
  clip width=4,
flip crossing/.list={3,5,8,13}
]
\strand([closed]0,1.5)..(-1.2,-1.5).. (3,-3.5) ..(0,-1.2) ..(-3,-3.5) ..(1.2,-1.5)..(0,1.5);
\strand ([closed]8,1.5)..(6.8,-1.5).. (11,-3.5) ..(8,-1.2) ..(5,-3.5) ..(9.2,-1.5)..(8,1.5);
\strand (1,-8)..(1.5,-8)..(2,-8) ..(2,-11)..(4,-12)..(6,-11)..(5,-10)..(6,-9)..(5,-8.5)..(3.5,-8.5)..(3,-9)..(4.5,-10.5)..(3.5,-11)..(3.3,-10.3)..(6,-10)..(5,-9)..(6,-7.5)..(4,-6.5)..(2,-7)..(1,-7); 

    \draw (3.7,-11.95) to (3.7,-11.1);
    \draw (3.8,-11.96) to (3.8,-11.1);
    \draw (1.3,-8) to (1.3,-6.95);
    \draw (1.4,-8) to (1.4,-6.95);
    \draw (5,-10.05) to (6.01,-10.05);
    \draw (5,-10.15) to (6.02,-10.15);
   
\end{knot}
\node[circle,draw=black, fill=black, inner sep=0pt,minimum size=5.5pt] (a) at (9.22,-1.33) {};
\end{tikzpicture}
\caption{Diagram of a classical knot, singular knot, and folded molecular chain.}
\label{diagrams}
\end{figure}
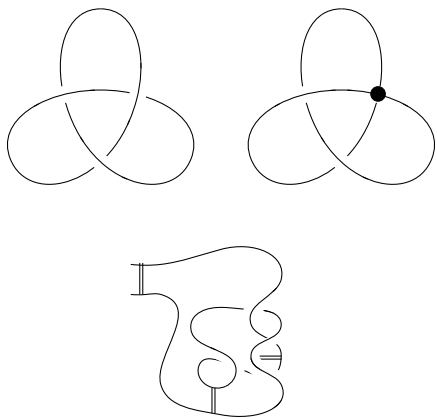

In Figure~\ref{diagrams}, we include a diagram of a classical knot, a singular knot and a folded molecular chain model.

We note that although the Reidemeister moves are essential in knot theory it is not practical to prove that two knot diagrams are not equivalent through the use of only Reidemeister moves. Instead we can compute and compare the values of a knot invariant . A \emph{knot invariant} is a map from the set $\mathcal{K}$ of knot diagrams to a fixed set $S$ whose values are constant on a knots' equivalence class. Suppose that $I: {\mathcal{K}} \rightarrow S$ is a knot invariant 
and let $K_1$ and $K_2$ be two diagrams.  If $K_1$ and $K_2$ are equivalent, then $I(K_1) = I(K_2)$. Thus, if $I(K_1) \neq I(K_2)$, then the knots $K_1$ and $K_2$ can not be equivalent.  The key idea is that a knot invariant will remain unchanged after applying Reidemeister moves to a knot diagram.  Notice that in some cases one may have $I(K_1) = I(K_2)$ without $K_1$ and $K_2$ being equivalent.  In this case one needs to search for a \emph{stronger} invariant (see Section~\ref{EBCI} below).    
In the case that $K_1$ and $K_2$ are equivalent if and only if $I(K_1) = I(K_2)$, then we have a \emph{complete invariant}.  Typical invariants are \emph{crossing number, unknotting number, Alexander and Jones polynomials}, etc. We will define specific invariants for knots, singular knots, and folded molecular chain models in the following section.

\section{Quandles, Singquandles and Bondles}\label{QSB}

In this section, we review algebraic structures useful for studying knots, singular knots and folded molecular chain models.  

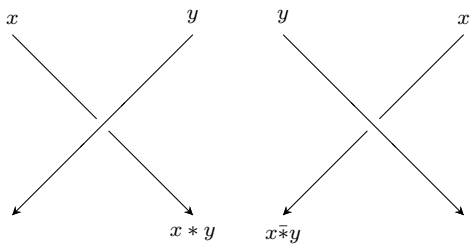
\begin{figure}[h]
\begin{tikzpicture}[use Hobby shortcut,scale=1.2]
\begin{knot}[
consider self intersections=true,
  ignore endpoint intersections=false,
  clip width=4,
  flip crossing/.list={3,4,5,6,11,13}
]
\strand[decoration={markings,mark=at position 1 with
    {\arrow[scale=1,>=stealth]{>}}},postaction={decorate}] (1,1) ..(-1,-1); 
\strand[decoration={markings,mark=at position 1 with
    {\arrow[scale=1,>=stealth]{>}}},postaction={decorate}] (-1,1) ..(1,-1);
    
\strand[decoration={markings,mark=at position 1 with
    {\arrow[scale=1,>=stealth]{>}}},postaction={decorate}] (2,1) ..(4,-1);
\strand[decoration={markings,mark=at position 1 with
    {\arrow[scale=1,>=stealth]{>}}},postaction={decorate}] (4,1) ..(2,-1);  

\end{knot}
\node[above] at (-1,1) {\tiny $x$};
\node[above] at (1,1) {\tiny $y$};
\node[below] at (1,-1) {\tiny $x*y$};
\node[above] at (2,1) {\tiny $y$};
\node[above] at (4,1) {\tiny $x$};
\node[below] at (2,-1) {\tiny $x\bar{*}y$};
\end{tikzpicture}
	\caption{Coloring rule at positive and negative crossings.}
		\label{Coloring}
\end{figure}
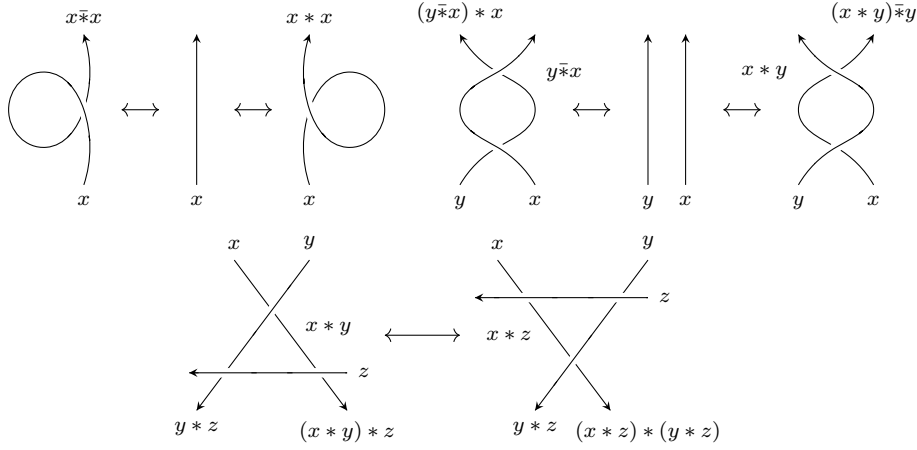
\begin{figure}[h]
\begin{tikzpicture}[use Hobby shortcut,add arrow/.style={postaction={decorate}, decoration={
  markings,
  mark=at position 1 with {\arrow[scale=1,>=stealth]{>}}}}]
\begin{knot}[
consider self intersections=true,
 ignore endpoint intersections=false,
  flip crossing/.list={1,5,6,7,8}
]

\strand(-6.5,-1)..(-6.5,0)..(-7,.5)..(-7.5,0)..(-7,-.5)..(-6.5,0)..(-6.5,1)[add arrow];
\draw[<->] (-6,0)..(-5.5,0);
\strand(-5,-1)..(-5,1)[add arrow];
\draw[<->] (-4.5,0)..(-4,0);
\strand(-3.5,-1)..(-3.5,0)..(-3,.5)..(-2.5,0)..(-3,-.5)..(-3.5,0)..(-3.5,1)[add arrow]; 
    

\strand(-.5,-1)..(-1,-.5)..(-1.5,0)..(-1,.5)..(-.5,1)[add arrow];   
\strand(-1.5,-1)..(-1,-.5)..(-.5,0)..(-1,.5)..(-1.5,1)[add arrow];
\draw[<->] (0,0).. (.5,0);
\strand(1,-1)..(1,1)[add arrow];    
\strand(1.5,-1)..(1.5,1)[add arrow];
\draw[<->] (2,0)..(2.5,0);
\strand(3,-1)..(3.5,-.5)..(4,0)..(3.5,.5)..(3,1)[add arrow];    
\strand(4,-1)..(3.5,-.5)..(3,0)..(3.5,.5)..(4,1)[add arrow];


\strand(-3,-3.5)..(-5.1,-3.5)[add arrow];
\strand(-3.5,-2)..(-5,-4)[add arrow];   
\strand(-4.5,-2)..(-3,-4)[add arrow];
 \draw[<->] (-2.5,-3)..(-1.5,-3);
\strand(1,-2.5)..(-1.3,-2.5)[add arrow];
\strand(1,-2)..(-.5,-4)[add arrow]; 
\strand(-1,-2)..(.5,-4)[add arrow];
\end{knot}
\node[below] at (-6.5,-1) {\tiny $x$};
\node[above] at (-6.5,1) {\tiny $x\bar{*}x$};
\node[below] at (-5,-1) {\tiny $x$};
\node[below] at (-3.5,-1) {\tiny $x$};
\node[above] at (-3.5,1) {\tiny $x*x$};
\node[below] at (1,-1) {\tiny $y$};
\node[below] at (1.5,-1) {\tiny $x$};
\node[below] at (-.5,-1) {\tiny $x$};
\node[below] at (-1.5,-1) {\tiny $y$};
\node[below] at (4,-1) {\tiny $x$};
\node[below] at (3,-1) {\tiny $y$};
\node[left] at (3,.5) {\tiny $x*y$};
\node[above] at (4,1) {\tiny $(x*y)\bar{*}y$};
\node[right] at (-.5,.5) {\tiny $y\bar{*}x$};
\node[above] at (-1.5,1) {\tiny $(y\bar{*}x)*x$};
\node[right] at (-3,-3.5) {\tiny $z$};
\node[above] at (-3.5,-2) {\tiny $y$};
\node[above] at (-4.5,-2) {\tiny $x$};
\node[right] at (-3.7,-2.9) {\tiny $x*y$};
\node[below] at (-3,-4) {\tiny $(x*y)*z$};
\node[below] at (-5,-4) {\tiny $y*z$};
\node[above] at (1,-2) {\tiny $y$};
\node[above] at (-1,-2) {\tiny $x$};
\node[right] at (1,-2.5) {\tiny $z$};
\node[left] at (-.4,-3) {\tiny $x*z$};
\node[below] at (-.5,-4) {\tiny $y*z$};
\node[below] at (1,-4) {\tiny $(x*z)*(y*z)$};
\end{tikzpicture}
	\caption{Coloring of Reidemeister moves.}
		\label{Axioms}
\end{figure}

    First, we start by reviewing the notion of colorings of knots by quandles (see \cite{EN} for more details).  A \emph{coloring} of a knot diagram corresponds to labeling each arc in the diagram with an element of a set $X$ subject to the coloring rules of Figure~\ref{Coloring}.  Next, by coloring the arcs of the three Reidemeister moves as in Figure~\ref{Axioms}, we obtain the necessary relations to define a quandle.

\begin{definition}
A set $X$ with \emph{two} binary operations, $*, \bar{*}:X \times X \rightarrow X$, satisfying the following $3$ axioms
\begin{enumerate}
\item
For all $x$ in $X$, $x*x=x$,
\item
For all $x,y$ in $X$, $(x*y)\bar{*}y=x=(x \bar{*}y)*y$,
    \item 
    For all $x,y,z$ in $X$,$ (x*y)*z=(x*z)*(y*z).$
\end{enumerate}
is called a \emph{quandle}.
\end{definition}

A \emph{typical} example of a quandle is the set $\mathbb{Z}_n$ of integers modulo $n$ with operations $x*y=-x+2y=x \bar{*}y$.  This example will be used later in Section~\ref{Examples} for some specific value(s) of the integer $n$.\\ 

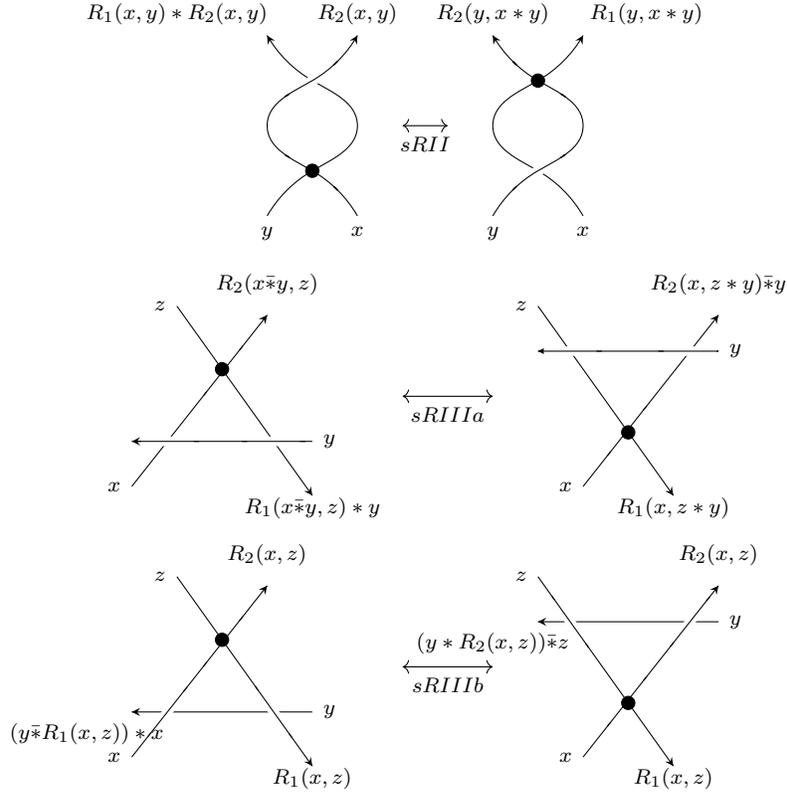
\begin{figure}[h]
\begin{tikzpicture}[use Hobby shortcut, scale=1.2, add arrow/.style={postaction={decorate}, decoration={
  markings,
  mark=at position 1 with {\arrow[scale=1,>=stealth]{>}}}}]
\begin{knot}[
consider self intersections=true,
 ignore endpoint intersections=false,
]
    

\strand(0,-1)..(-.5,-.5)..(-1,0)..(-.5,.5)..(0,1)[add arrow];   
\strand(-1,-1)..(-.5,-.5)..(0,0)..(-.5,.5)..(-1,1)[add arrow];
\draw[<->] (.5,0).. (1,0);

\strand(1.5,-1)..(2,-.5)..(2.5,0)..(2,.5)..(1.5,1)[add arrow];    
\strand(2.5,-1)..(2,-.5)..(1.5,0)..(2,.5)..(2.5,1)[add arrow];


\strand(-.5,-3.5)..(-2.5,-3.5)[add arrow];
\strand(-2.5,-4)..(-1,-2.1)[add arrow];   
\strand(-2,-2)..(-.5,-4.1)[add arrow];
 \draw[<->] (.5,-3)..(1.5,-3);
\strand(4,-2.5)..(2,-2.5)[add arrow];
\strand(2.5,-4)..(4,-2.1)[add arrow]; 
\strand(2,-2)..(3.5,-4.1)[add arrow];


\strand(-2.5,-7)..(-1,-5.1)[add arrow];   
\strand(-2,-5)..(-.5,-7.1)[add arrow];
\strand(-.5,-6.5)..(-2.5,-6.5)[add arrow];
 \draw[<->] (.5,-6)..(1.5,-6);
\strand(2.5,-7)..(4,-5.1)[add arrow]; 
\strand(2,-5)..(3.5,-7.1)[add arrow];
\strand(4,-5.5)..(2,-5.5)[add arrow];
\end{knot}


\node[below] at (0,-1) {\tiny $x$};
\node[below] at (-1,-1) {\tiny $y$};
\node[above] at (0,1) {\tiny $R_2(x,y)$};
\node[above] at (-2,1) {\tiny $R_1(x,y)*R_2(x,y)$};

\node[below] at (2.5,-1) {\tiny $x$};
\node[below] at (1.5,-1) {\tiny $y$};
\node[above] at (1.5,1) {\tiny $R_2(y,x*y)$};
\node[above] at (3.2,1) {\tiny $R_1(y,x*y)$};

\node[right] at (-.5,-3.5) {\tiny $y$};
\node[left] at (-2.5,-4) {\tiny $x$};
\node[left] at (-2,-2) {\tiny $z$};
\node[below] at (-.5,-4) {\tiny $R_1(x\bar{*}y,z)*y$};
\node[above] at (-1,-2) {\tiny $R_2(x\bar{*}y,z)$};

\node[right] at (4,-2.5) {\tiny $y$};
\node[left] at (2.5,-4) {\tiny $x$};
\node[left] at (2,-2) {\tiny $z$};
\node[below] at (3.5,-4) {\tiny $R_1(x,z*y)$};
\node[above] at (4,-2) {\tiny $R_2(x,z*y)\bar{*}y$};

\node[right] at (-.5,-6.5) {\tiny $y$};
\node[left] at (-2.5,-7) {\tiny $x$};
\node[left] at (-2,-5) {\tiny $z$};
\node[below] at (-.5,-7) {\tiny $R_1(x,z)$};
\node[above] at (-1,-5) {\tiny $R_2(x,z)$};
\node[below] at (-3,-6.5) {\tiny $(y\bar{*}R_1(x,z))*x$};

\node[right] at (4,-5.5) {\tiny $y$};
\node[left] at (2.5,-7) {\tiny $x$};
\node[left] at (2,-5) {\tiny $z$};
\node[below] at (3.5,-7) {\tiny $R_1(x,z)$};
\node[above] at (4,-5) {\tiny $R_2(x,z)$};
\node[below] at (1.5,-5.5) {\tiny $(y*R_2(x,z))\bar{*}z$};

\node[below] at (.75,0) {\tiny $sRII$};
\node[below] at (1,-3) {\tiny $sRIIIa$};
\node[below] at (1,-6) {\tiny $sRIIIb$};

\node[circle,draw=black, fill=black, inner sep=0pt,minimum size=5pt] (a) at (-.5,-.5) {};
\node[circle,draw=black, fill=black, inner sep=0pt,minimum size=5pt] (a) at (2,.5) {};

\node[circle,draw=black, fill=black, inner sep=0pt,minimum size=5pt] (a) at (-1.5,-2.7) {};
\node[circle,draw=black, fill=black, inner sep=0pt,minimum size=5pt] (a) at (3,-3.4) {};

\node[circle,draw=black, fill=black, inner sep=0pt,minimum size=5pt] (a) at (-1.5,-5.7) {};
\node[circle,draw=black, fill=black, inner sep=0pt,minimum size=5pt] (a) at (3,-6.4) {};
\end{tikzpicture}
	\caption{Additional singular Reidemeister moves with colorings, $sRII$, $sRIIIa$ and $sRIIIb$.}
		\label{singmoves}
\end{figure}

 In \cite{CEHN}, the authors introduced an algebraic structure to study singular knots called an \emph{oriented singquandle}. Similar to Reidemeister moves in classical knot theory, it is suitable to consider the generalized Reidemeister moves \cite{BEHY} in order to obtain an invariant of singular knots. We recall the definition of singquandle whose axioms come from the classical Reidemeister moves and the additional singular Reidemeister moves, see Figure~\ref{singmoves}.
 
\begin{definition}\label{oriented SingQdle}
	Let $(X, *)$ be a quandle.  Let $R_1$ and $R_2$ be two maps from $X \times X$ to $X$.  The quadruple $(X, *, R_1, R_2)$ is called an {\it oriented singquandle} if the following axioms are satisfied for all $x,y,z \in X$:
	\begin{eqnarray}
		R_1(x \bar{*}y,z)*y&=&R_1(x,z*y)  \label{eq1}\\
		R_2(x\bar{*}y, z) & =&  R_2(x,z*y)\bar{*}y \label{eq2}\\
	      (y\bar{*}R_1(x,z))*x   &=& (y*R_2(x,z))\bar{*}z  \label{eq3}\\
R_2(x,y)&=&R_1(y,x*y)   \label{eq4}\\
R_1(x,y)*R_2(x,y)&=&R_2(y,x*y).   \label{eq5}	
\end{eqnarray}	
\end{definition}

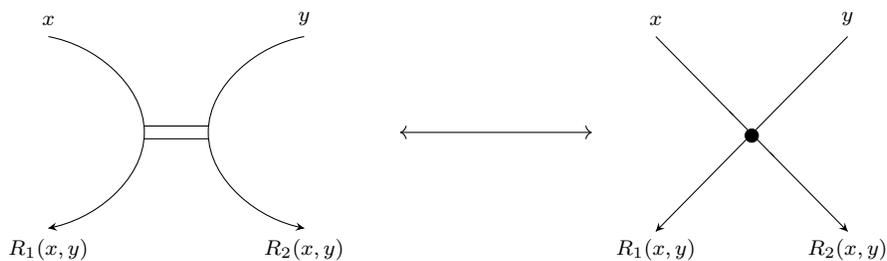
\begin{figure}[h]
\begin{tikzpicture}[use Hobby shortcut,scale=.85,add arrow/.style={postaction={decorate}, decoration={
  markings,
  mark=at position 1 with {\arrow[scale=1,>=stealth]{>}}}}]
\begin{knot}

 \strand (-3.5,1.5)..(-2.5,1)..(-2,0)..(-2.5,-1)..(-3.5,-1.5)[add arrow];
 \strand (.5,1.5)..(-.5,1)..(-1,0)..(-.5,-1) ..(.5,-1.5)[add arrow];

\draw[<->] (2,0)..(5,0);

\strand (6,1.5) ..(9,-1.55) [add arrow];
\strand (9,1.5) ..(6,-1.55) [add arrow]; 
\end{knot}

\node[circle,draw=black, fill=black, inner sep=0pt,minimum size=5pt] (a) at (7.5,-.05) {};

 \draw (-2,.1) to (-1,.1);
 \draw (-2,-.1) to (-1,-.1);

\node[above] at (-3.5,1.5) {\tiny $x$};
\node[above] at (.5,1.5) {\tiny $y$};
\node[below] at (-3.5,-1.5) {\tiny $R_1(x,y)$};
\node[below] at (.5,-1.5) {\tiny $R_2(x,y)$};

\node[above] at (6,1.5) {\tiny $x$};
\node[above] at (9,1.5) {\tiny $y$};
\node[below] at (6,-1.5) {\tiny $R_1(x,y)$};
\node[below] at (9,-1.5) {\tiny $R_2(x,y)$};
\end{tikzpicture}
\vspace{.2in}
		\caption{Equivalence of a bond with parallel strands in a folded molecular chain diagram and a singular crossing in a singular knot diagram.}
		\label{bondsing}
\end{figure}

\begin{figure}[h]
\begin{tikzpicture}[use Hobby shortcut,scale=.85,add arrow/.style={postaction={decorate}, decoration={
  markings,
  mark=at position 1 with {\arrow[scale=1,>=stealth]{>}}}}]
\begin{knot}
 \strand (3,1.5)..(4,1)..(4.5,0)..(4,-1)..(3,-1.5)[add arrow];
 \strand (7,-1.5)..(6,-1)..(5.5,0)..(6,1)..(7,1.5)[add arrow];

\end{knot}

 \draw (4.5,.1) to (5.5,.1);
 \draw (4.5,-.1) to (5.5,-.1);
\node[above] at (3,1.5) {\tiny $x$};
\node[below] at (7,-1.5) {\tiny $y$};
\node[below] at (3,-1.5) {\tiny $R_3(x,y)$};
\node[above] at (7,1.5) {\tiny $R_3(y,x)$};

\end{tikzpicture}
\vspace{.2in}
		\caption{Bond with anti-parallel strands.}
		\label{bondanti}
\end{figure}
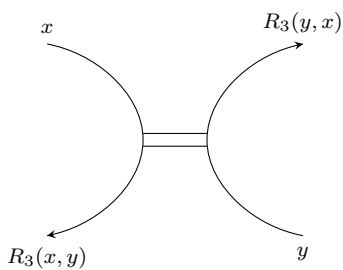

Note that oriented singquandles are useful for classifying singular knots. In this paper, we are interested in distinguishing folded molecular chains. Therefore, we introduce an extension of oriented singquandles defined in \cite{ADEM}. The authors extended the definition of an oriented singquandle to take into account the bonds in a folded molecular chain. To this end, the authors used singular crossings to represent bonds with parallel strands, see Figure~\ref{bondsing}. Furthermore, the authors assigned a new function to bonds with anti-parallel strands, see Figure~\ref{bondanti}. The inclusion of the bonds with anti-parallel strands to the folded molecular diagram requires three more Reidemeister moves. We will refer to them as the \emph{bond Reidemeister moves}. The generalization of an oriented singquandle must therefore satisfy the four additional relations obtained from the bond Reidemeister moves, see Figure~\ref{bondle}.

\begin{definition}
An \emph{oriented bondle} is an oriented singquandle $(X,*, R_1, R_2)$ and a choice for a function $R_3(x,y)$ such that the following relations are satisfied for $x,y,z \in X$:

\begin{eqnarray}
R_3(y,x\bar{*} z) &=& R_3(y*z,x)\bar{*}z\\
R_3(x,y*z) &=& R_3(x\bar{*}z,y)*z\\
(z \bar{*} R_3(x,y))*x &=& (z \bar{*}y)*R_3(y,x)\\
R_3(x,y)\bar{*} y &=& R_3(x\bar{*}R_3(y,x),y).
\end{eqnarray}
\end{definition}
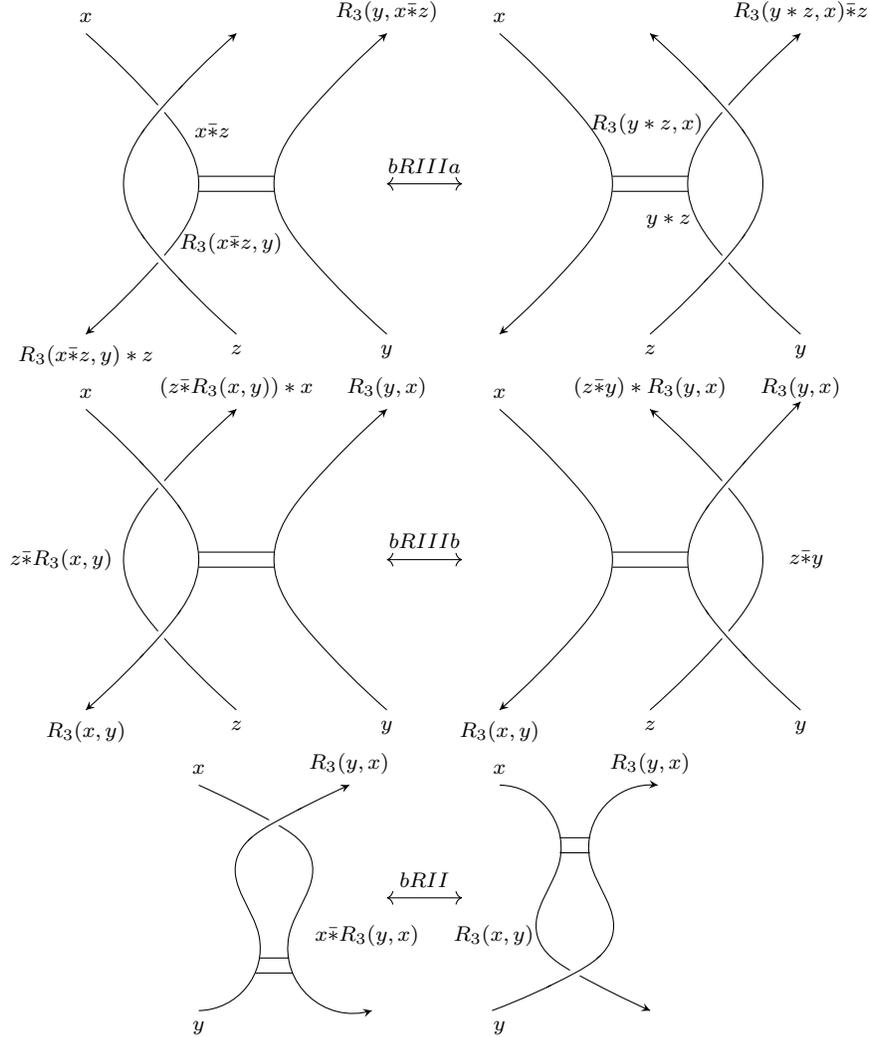
\begin{figure}[h]
\begin{tikzpicture}[use Hobby shortcut,scale=1,add arrow/.style={postaction={decorate}, decoration={
  markings,
  mark=at position 1 with {\arrow[scale=1,>=stealth]{>}}}}]
\begin{knot}
 \strand (-1.5,-2)..(-2.5,-1)..(-3,0)..(-2.5,1)..(-1.5,2)[add arrow];
 \strand (-3.5,2)..(-2.5,1)..(-2,0)..(-2.5,-1)..(-3.5,-2)[add arrow];
 \strand (.5,-2)..(-.5,-1)..(-1,0)..(-.5,1)..(.5,2)[add arrow];
 
 \draw[<->] (.5,0)..(1.5,0);
 \strand (4,-2)..(5,-1)..(5.5,0)..(5,1)..(4,2)[add arrow];
 \strand (2,2)..(3,1)..(3.5,0)..(3,-1)..(2,-2)[add arrow];
 \strand (6,-2)..(5,-1)..(4.5,0)..(5,1)..(6,2)[add arrow];
 
 \draw (-2,.1) to (-1,.1);
 \draw (-2,-.1) to (-1,-.1);
 
 \draw (3.5,.1) to (4.5,.1);
 \draw (3.5,-.1) to (4.5,-.1);
 
  
 \strand (-3.5,-3)..(-2.5,-4)..(-2,-5)..(-2.5,-6)..(-3.5,-7)[add arrow];
 \strand (.5,-7)..(-.5,-6)..(-1,-5)..(-.5,-4)..(.5,-3)[add arrow];
 \strand (-1.5,-7)..(-2.5,-6)..(-3,-5)..(-2.5,-4)..(-1.5,-3)[add arrow];

 \draw[<->] (.5,-5)..(1.5,-5);

 \strand (6,-7)..(5,-6)..(4.5,-5)..(5,-4)..(6,-2.9)[add arrow]; 
 \strand (4,-7)..(5,-6)..(5.5,-5)..(5,-4)..(4,-3)[add arrow];
 \strand (2,-3)..(3,-4)..(3.5,-5)..(3,-6)..(2,-7)[add arrow];
 
 \draw (-2,-4.9) to (-1,-4.9);
 \draw (-2,-5.1) to (-1,-5.1);
 
 \draw (3.5,-4.9) to (4.5,-4.9);
 \draw (3.5,-5.1) to (4.5,-5.1);
 
 \strand (-2,-11).. (-1.2,-10)..(-1.5,-9)..(-1,-8.5)..(0,-8)[add arrow];
 \strand (-2,-8)..(-1,-8.5)..(-.5,-9)..(-.8,-10)..(0.3,-11)[add arrow];

 \draw (-1.2,-10.3) to (-.8,-10.3);
 \draw (-1.25,-10.5) to (-.75,-10.5);
 
 \draw[<->] (.5,-9.5)..(1.5,-9.5);
 
 \strand (1.9,-11)..(3,-10.5)..(3.5,-10)..(3.2,-9)..(4.1,-8)[add arrow];
 \strand (2,-8)..(2.8,-9)..(2.5,-10)..(3,-10.5)..(4,-11)[add arrow];

 \draw (2.8,-8.9) to (3.2,-8.9);
 \draw (2.8,-8.7) to (3.2,-8.7);

\end{knot}
\node[below] at (-1.5,-2) {\tiny $z$};
\node[above] at (-3.5,2) {\tiny $x$};
\node[below] at (.5,-2) {\tiny $y$};
\node[right] at (-2.2,.7) {\tiny $x\bar{*}z$};
\node[right] at (-2.4,-.8) {\tiny $R_3(x \bar{*}z,y)$};
\node[above] at (.5,2) {\tiny $R_3(y,x\bar{*}z)$};
\node[below] at (-3.5,-2) {\tiny $R_3(x\bar{*}z,y)*z$};

\node[below] at (4,-2) {\tiny $z$};
\node[above] at (2,2) {\tiny $x$};
\node[below] at (6,-2) {\tiny $y$};
\node[right] at (3.8,-.5) {\tiny $y*z$};
\node[left] at (4.85,.79) {\tiny $R_3(y*z,x)$};
\node[above] at (6,2) {\tiny $R_3(y*z,x)\bar{*}z$};

\node[below] at (-1.5,-7) {\tiny $z$};
\node[above] at (-3.5,-3) {\tiny $x$};
\node[below] at (.5,-7) {\tiny $y$};
\node[left] at (-3,-5) {\tiny $z\bar{*}R_3(x,y)$};
\node[below] at (-3.5,-7) {\tiny $R_3(x,y)$};
\node[above] at (.5,-3) {\tiny $R_3(y,x)$};
\node[above] at (-1.5,-3) {\tiny $(z\bar{*}R_3(x,y))*x$};

\node[below] at (4,-7) {\tiny $z$};
\node[above] at (2,-3) {\tiny $x$};
\node[below] at (6,-7) {\tiny $y$};
\node[right] at (5.7,-5) {\tiny $z\bar{*}y$};
\node[above] at (4,-3) {\tiny $(z\bar{*}y)*R_3(y,x)$};
\node[below] at (2,-7) {\tiny $R_3(x,y)$};
\node[above] at (6,-3) {\tiny $R_3(y,x)$};

\node[above] at (-2,-8) {\tiny $x$};
\node[below] at (-2,-11) {\tiny $y$};
\node[above] at (0,-8) {\tiny $R_3(y,x)$};
\node[right] at (-.6,-10) {\tiny $x \bar{*} R_3(y,x)$};

\node[above] at (2,-8) {\tiny $x$};
\node[below] at (2,-11) {\tiny $y$};
\node[above] at (4,-8) {\tiny $R_3(y,x)$};
\node[left] at (2.6,-10) {\tiny $R_3(x,y)$};

\node[above] at (1,0) {\tiny $bRIIIa$};
\node[above] at (1,-5) {\tiny $bRIIIb$};
\node[above] at (1,-9.5) {\tiny $bRII$};

\end{tikzpicture}

\caption{Oriented bond Reidemeister moves with bondle colorings, $bRIIIa$, $bRIIIb$, and $bRII$.}
\label{bondle}
\end{figure}


In \cite{ADEM}, it was shown that oriented bondles give an invariant of folded molecular chains.  Precisely, it was shown that given two projections of folded molecular chains and a fixed bondle, if the numbers of distinct colorings of each projection by the bondle are distinct then the two folded molecular chains are not topologically equivalent.

The following result was proved in \cite{ADEM}.
\begin{proposition}\label{bondleEx}
Let $n =pq$ where $p$ and $q$ are odd primes. Assume further that $x*y= ax + (1-a)y$ and $x\bar{*}y=a^{-1}x+(1-a^{-1})y$, for invertible element $a$ in $\mathbb{Z}_n$. For any fixed element $b$ in $\mathbb{Z}_n$, let $R_1(x,y)=bx +(1-b)y$, $R_2(x,y)= a(1-b)x+\left[ b+(1-a)(1-b) \right]y$ and $R_3(x,y)=mx + (1-m)y$. Then $(\mathbb{Z}_n,*, R_1, R_2, R_3)$ is an oriented bondle if and only if $p$ divides $m$ and $q$ divides $(m-1)$ or $p$ divides $(m-1)$ and $q$ divides $m$.
\end{proposition}
Note that Proposition~\ref{bondleEx} defines a nice family of bondles. Specifically, the above proposition defines bondles with linear operations $*, R_1,R_2$ and $R_3$. Having linear operations will simplify computations. In the following example we will define a specific bondle using Proposition~\ref{bondleEx} and then compute the bondle counting invariant for a specific folded molecular chain.

\begin{example}
Let $(\mathbb{Z}_{15}, *, R_1,R_2,R_3)$ be an oriented bondle with $x*y = 8(x+y)$, $R_1(x,y) = 2x - y$, and $R_2(x,y) = 7x-6y$. Note that $R_3(x,y)$ need not be defined since our diagram does not have any anti-parallel bonds. 
\begin{figure}[h]
\begin{tikzpicture}[use Hobby shortcut,scale=.7]
\begin{knot}[
consider self intersections=true,
clip width=4,
  ignore endpoint intersections=false,
]
\strand[decoration={markings,mark=at position .1 with
    {\arrow[scale=1.5,>=stealth]{>}}},postaction={decorate}] (-1,1)..(-.5,1)..(0,1)..(1,-2)..(2,-1)..(2,1)..(3,3)..(4,3)..(5,2)..(5,0)..(5,-1)..(4,-2)..(3,-1)..(3,1)..(2,3).. (1,3)..(-.5,2)..(-1,2); 
    \draw (-.7,1) to (-.7,1.97);
    \draw (-.6,1) to (-.6,2);
    \draw (2,-1.05) to (3,-1.05);
    \draw (1.96,-1.15) to (3.05,-1.15);
\end{knot}

\node[below] at (1,-2) {\tiny $x$};
\node[below] at (4,-2) {\tiny $y$};
\node at (1,1.5) {\tiny $R_2(y,x)$};
\node at (4,1.5) {\tiny $R_1(y,x)$};
\end{tikzpicture}
\caption{Coloring of $P$}
\end{figure}
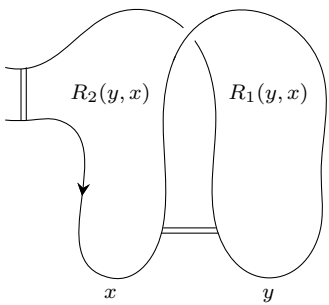
A coloring of $P$ gives the following equation
\[ R_2(y,x) =  y \]
which simplifies
\[ 6(y-x) = 0. \]
So if $5$ divides $y-x$, that is, $y=x$ (trivial colorings) or $y=x+5$ or $y=x+10$ (these last 2 cases correspond to nontrivial colorings).   Thus, the total number of solutions is $3 \times 15=45$. This means that this protein has 45 colorings by our chosen bondle. 
\end{example}
\section{Enhancement of the Bondle Counting Invariant}\label{EBCI}

In this section, we introduce a function at bonds with anti-parallel strands, which we will call a \emph{Boltzmann weight}. We will also introduce previously defined functions for classical crossings and singular crossings (bonds with parallel strands). We will use the Boltzmann weights to enhance the bondle coloring invariant. 
We will first review the construction of quandle  2-cocycles (see \cite{EN} page 192) and the Boltzmann weights of singquandles (see \cite{CCEH} for more details).

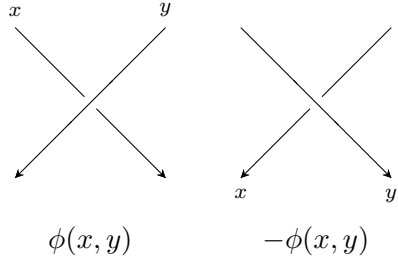
\begin{figure}[h]
\begin{tikzpicture}[use Hobby shortcut]
\begin{knot}[
consider self intersections=true,
clip width=4,
  ignore endpoint intersections=false,
  flip crossing/.list={3,4,5,6,11,13}
]
\strand[decoration={markings,mark=at position 1 with
    {\arrow[scale=1,>=stealth]{>}}},postaction={decorate}] (1,1) ..(-1,-1); 
\strand[decoration={markings,mark=at position 1 with
    {\arrow[scale=1,>=stealth]{>}}},postaction={decorate}] (-1,1) ..(1,-1);
    
\strand[decoration={markings,mark=at position 1 with
    {\arrow[scale=1,>=stealth]{>}}},postaction={decorate}] (2,1) ..(4,-1);
\strand[decoration={markings,mark=at position 1 with
    {\arrow[scale=1,>=stealth]{>}}},postaction={decorate}] (4,1) ..(2,-1);  

\end{knot}
\node[above] at (-1,1) {\tiny $x$};
\node[above] at (1,1) {\tiny $y$};
\node[below] at (0,-1.5) {\small $\phi(x,y)$};
\node[below] at (2,-1) {\tiny $x$};
\node[below] at (4,-1) {\tiny $y$};
\node[below] at (3,-1.5) {\small $-\phi(x,y)$};
\end{tikzpicture}
	\caption{Weight functions at a positive and negative crossing.}
		\label{weights}
\end{figure}
Let $(X,*)$ be a quandle, $A$ be an abelian group and $\phi$ be a function from $X \times X$ to $A$.  For $x,y \in X$, we will assign the function values, $\phi(x,y)$ at a positive crossing and $-\phi(x,y)$ at a negative crossing, as shown in Figure~\ref{weights}. The function $\phi: X \times X \rightarrow A$ is called a \emph{quandle 2-cocyle} if it satisfies the following condition coming from Reidemeister move III.  
\begin{equation}\label{cocycle1}
 \phi(x,y) + \phi(x*y,z) = \phi(x,z) + \phi(x*z,y*z), \quad \text{for all $x,y \in X$}. 
\end{equation}
Additionally, Reidemeister move RI imposes the following condition:
\begin{equation}\label{cocycle2}
\phi(x,x) = 0, \quad \text{for all $x\in X$}.
\end{equation}
Lastly, Reidemeister move RII is automatically satisfied since the weights at positive and negative crossing cancel each other.

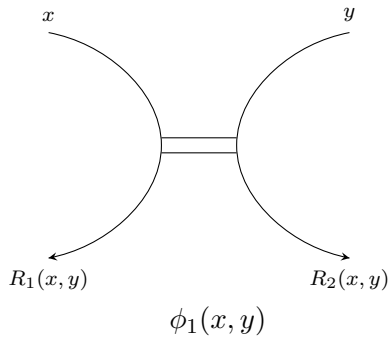
\begin{figure}[h]
\begin{tikzpicture}[use Hobby shortcut,add arrow/.style={postaction={decorate}, decoration={
  markings,
  mark=at position 1 with {\arrow[scale=1,>=stealth]{>}}}}]
\begin{knot}

 \strand (-3.5,1.5)..(-2.5,1)..(-2,0)..(-2.5,-1)..(-3.5,-1.5)[add arrow];
 \strand (.5,1.5)..(-.5,1)..(-1,0)..(-.5,-1) ..(.5,-1.5)[add arrow];

\end{knot}

 \draw (-2,.1) to (-1,.1);
 \draw (-2,-.1) to (-1,-.1);

\node[above] at (-3.5,1.5) {\tiny $x$};
\node[above] at (.5,1.5) {\tiny $y$};
\node[below] at (-3.5,-1.5) {\tiny $R_1(x,y)$};
\node[below] at (.5,-1.5) {\tiny $R_2(x,y)$};
\node[below] at (-1.25,-2) {\small $\phi_1(x,y)$};

\end{tikzpicture}
\vspace{.2in}
		\caption{Boltzmann weight at a bond/singular crossing.}
		\label{singweight}
\end{figure}

In \cite{CCEH}, the 2-cocycle idea was extended to singular knots and links. The authors introduced a second function $\phi_1$ and assigned this function to singular crossings as in Figure~\ref{singweight}. We will adopt the same convention, but recall that singular crossings represent bonds with parallel strands in projections of folded molecular chain models. 

Let $(X,*,R_1,R_2)$ be an oriented singquandle and let $\phi$ be a 2-cocycle with values in $A$. We will will consider a new function $\phi_1: X \times X \rightarrow A$ to assign a Boltzmann weight at a singular crossing. The functions $\phi$ and $\phi_1$ will satisfy some conditions imposed by the additional singular Reidemeister moves. Therefore, Reidemeister move $sRIIIa$ in Figure~\ref{singmoves} imposes the following condition on $\phi$ and $\phi_1$, for all $x,y,z \in X$:
\begin{equation}\label{sing1}
\begin{split}
-\phi(x\bar{*}y,y)+\phi_1(x\bar{*}y,z)+\phi(R_1(x\bar{*}y,z),y)=\\\phi(z,y)+\phi_1(x,z*y)-\phi(R_2(x,z*y)\bar{*}y,y).
\end{split}
\end{equation}
From Reidemeister move $sRIIIb$ in Figure~\ref{singmoves} we obtain the following condition, for all $x,y,z \in X$:
\begin{equation}\label{sing2}
\begin{split}
  \phi(y\bar{*}R_1(x,z),x)-\phi(y\bar{*}R_1(x,z),R_1(x,z))=\\
-\phi((y*R_2(x,z))\bar{*}z, z)+\phi(y,R_2(x,z)). 
\end{split}
\end{equation}
Lastly, Reidemeister move $sRII$ in Figure~\ref{singmoves} imposes the following condition, for all $x,y \in X$,

\begin{equation}\label{sing3}
\begin{split}
\phi_1(x,y)+\phi(R_1(x,y),R_2(x,y))=
\phi(x,y)+\phi_1(y,x*y). 
\end{split}
\end{equation}

\begin{figure}[h]
\begin{tikzpicture}[use Hobby shortcut,add arrow/.style={postaction={decorate}, decoration={
  markings,
  mark=at position 1 with {\arrow[scale=1,>=stealth]{>}}}}]
\begin{knot}
 \strand (3,1.5)..(4,1)..(4.5,0)..(4,-1)..(3,-1.5)[add arrow];
 \strand (7,-1.5)..(6,-1)..(5.5,0)..(6,1)..(7,1.5)[add arrow];

 \strand (9,-1.5)..(10,-1)..(10.5,0)..(10,1)..(9,1.5)[add arrow];
 \strand (13,1.5)..(12,1)..(11.5,0)..(12,-1)..(13,-1.5)[add arrow];
 
\end{knot}

 \draw (4.5,.1) to (5.5,.1);
 \draw (4.5,-.1) to (5.5,-.1);

 \draw (10.5,.1) to (11.5,.1);
 \draw (10.5,-.1) to (11.5,-.1);

\node[above] at (3,1.5) {\tiny $x$};
\node[below] at (7,-1.5) {\tiny $y$};
\node[below] at (3,-1.5) {\tiny $R_3(x,y)$};
\node[above] at (7,1.5) {\tiny $R_3(y,x)$};
\node[below] at (5,-2) {\small $\phi_2(x,y)$};

\node[below] at (9,-1.5) {\tiny $x$};
\node[above] at (13,1.5) {\tiny $y$};
\node[above] at (9,1.5) {\tiny $R_3(x,y)$};
\node[below] at (13,-1.5) {\tiny $R_3(y,x)$};
\node[below] at (11,-2) {\small $\phi_2(x,y)$};
\end{tikzpicture}
\vspace{.2in}
		\caption{Boltzmann weight at a bond with anti-parallel strands.}
		\label{bondweight}
\end{figure}
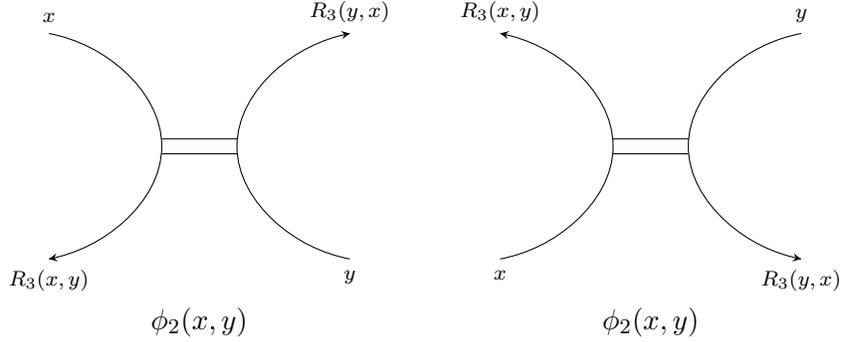

Note that if $(X,*,R_1,R_2)$ is an oriented singquandle, and $\phi, \phi_1 : X\times X \rightarrow A$ satisfy equations (10)-(14) we will refer to $\phi$ and $\phi_1$ as Boltzmann weights. In \cite{CCEH}, the Boltzmann weights of an oriented singquandle were used to define a proper enhancement of the singquandle coloring invariant. Since we are interested in folded molecular chain models, we will define Boltzmann weights for oriented bondles. Let $(X,*,R_1,R_2,R_3)$ be an oriented bondle with Boltzmann weight $\phi,\phi_1: X \times X \rightarrow A$ with $\phi$ assigned to classical crossings and $\phi_1$ assigned to bonds with parallel strands. We now define a third Boltzmann weight to allow contributions from bonds with anti-parallel strands. Let $\phi_2 : X \times X \rightarrow A$, we will follow the convention in Figure~\ref{bondweight}. 

\begin{figure}[h]
\begin{tikzpicture}[use Hobby shortcut,add arrow/.style={postaction={decorate}, decoration={
  markings,
  mark=at position 1 with {\arrow[scale=1,>=stealth]{>}}}}]
\begin{knot}
 \strand (-1.5,-2)..(-2.5,-1)..(-3,0)..(-2.5,1)..(-1.5,2)[add arrow];
 \strand (-3.5,2)..(-2.5,1)..(-2,0)..(-2.5,-1)..(-3.5,-2)[add arrow];
 \strand (.5,-2)..(-.5,-1)..(-1,0)..(-.5,1)..(.5,2)[add arrow];
 
 \draw[<->] (.5,0)..(1.5,0);
 \strand (4,-2)..(5,-1)..(5.5,0)..(5,1)..(4,2)[add arrow];
 \strand (2,2)..(3,1)..(3.5,0)..(3,-1)..(2,-2)[add arrow];
 \strand (6,-2)..(5,-1)..(4.5,0)..(5,1)..(6,2)[add arrow];
 
 \draw (-2,.1) to (-1,.1);
 \draw (-2,-.1) to (-1,-.1);
 
 \draw (3.5,.1) to (4.5,.1);
 \draw (3.5,-.1) to (4.5,-.1);

\end{knot}

\node[below,gray] at (-1.5,-2) {\tiny $z$};
\node[above,gray] at (-3.5,2) {\tiny $x$};
\node[below,gray] at (.5,-2) {\tiny $y$};
\node[right,gray] at (-2.2,.7) {\tiny $x\bar{*}z$};
\node[right,gray] at (-2.4,-.8) {\tiny $R_3(x \bar{*}z,y)$};
\node[above,gray] at (.5,2) {\tiny $R_3(y,x\bar{*}z)$};
\node[below,gray] at (-3.5,-2) {\tiny $R_3(x\bar{*}z,y)*z$};
\node[left] at (-2.55,1) {\tiny $-\phi(x\bar{*}z,z)$};
\node[left] at (-2.55,-1) {\tiny $\phi(R_3(x\bar{*}z,y),z)$};
\node[right] at (-.9,.3) {\tiny $\phi_2(x\bar{*}z,y)$};
\node[below,gray] at (4,-2) {\tiny $z$};
\node[above,gray] at (2,2) {\tiny $x$};
\node[below,gray] at (6,-2) {\tiny $y$};
\node[right,gray] at (3.8,-.5) {\tiny $y*z$};
\node[left,gray] at (4.85,.79) {\tiny $R_3(y*z,x)$};
\node[above,gray] at (6,2) {\tiny $R_3(y*z,x)\bar{*}z$};
\node[right] at (5.3,1) {\tiny $-\phi(R_3(y*z,x)\bar{*}z,z)$};
\node[right] at (5.5,-1) {\tiny $\phi(y,z)$};
\node[left] at (3.3,.3) {\tiny $\phi_2(x,y*z)$};
\end{tikzpicture}
\vspace{.2in}
		\caption{Bondle coloring and Boltzmann weight for $bRIIIa$.}
		\label{weightax1}
\end{figure}
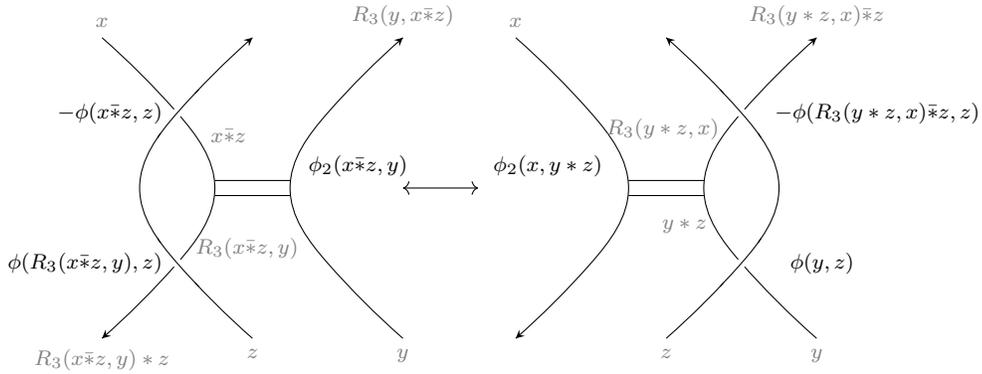

We require that $\phi_2$ is compatible with $\phi$ and $\phi_1$, therefore, the weight functions must satisfy the relations imposed by bond Reidemeister moves. Using Figure~\ref{weightax1}, we obtain the following condition:
\begin{equation}\label{bondle1}
\begin{split}
 -\phi(x \bar{*} z, z) + \phi(R_3(x\bar{*}z,y),z) + \phi_2(x\bar{*}z,y) =\\ \phi_2(x, y * z) - \phi(R_3(y*z, x) \bar{*}z,z)+ \phi(y,z)\quad \text{for all $x,y,z \in X$}.     
\end{split}
\end{equation}

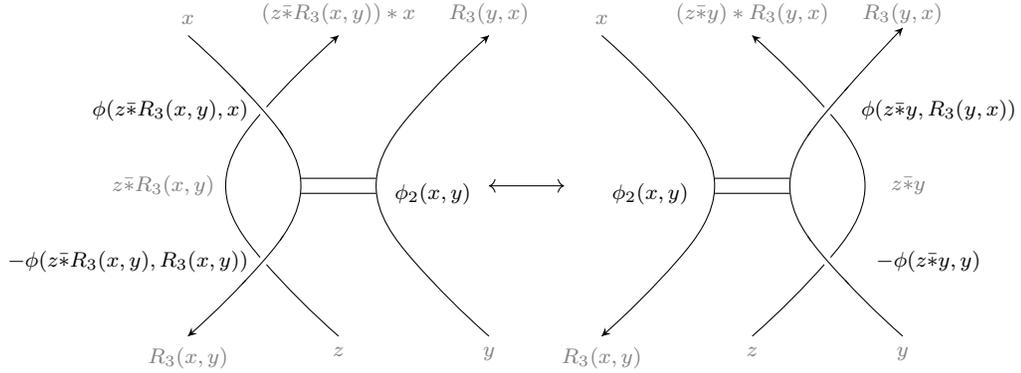
\begin{figure}[h]
\begin{tikzpicture}[use Hobby shortcut,add arrow/.style={postaction={decorate}, decoration={
  markings,
  mark=at position 1 with {\arrow[scale=1,>=stealth]{>}}}}]
\begin{knot}

  
 \strand (-3.5,-3)..(-2.5,-4)..(-2,-5)..(-2.5,-6)..(-3.5,-7)[add arrow];
 \strand (.5,-7)..(-.5,-6)..(-1,-5)..(-.5,-4)..(.5,-3)[add arrow];
 \strand (-1.5,-7)..(-2.5,-6)..(-3,-5)..(-2.5,-4)..(-1.5,-3)[add arrow];

 \draw[<->] (.5,-5)..(1.5,-5);

 \strand (6,-7)..(5,-6)..(4.5,-5)..(5,-4)..(6,-2.9)[add arrow]; 
 \strand (4,-7)..(5,-6)..(5.5,-5)..(5,-4)..(4,-3)[add arrow];
 \strand (2,-3)..(3,-4)..(3.5,-5)..(3,-6)..(2,-7)[add arrow];
 
 \draw (-2,-4.9) to (-1,-4.9);
 \draw (-2,-5.1) to (-1,-5.1);
 
 \draw (3.5,-4.9) to (4.5,-4.9);
 \draw (3.5,-5.1) to (4.5,-5.1);

\end{knot}
\node[below,gray] at (-1.5,-7) {\tiny $z$};
\node[above,gray] at (-3.5,-3) {\tiny $x$};
\node[below,gray] at (.5,-7) {\tiny $y$};
\node[left,gray] at (-3,-5) {\tiny $z\bar{*}R_3(x,y)$};
\node[below,gray] at (-3.5,-7) {\tiny $R_3(x,y)$};
\node[above,gray] at (.5,-3) {\tiny $R_3(y,x)$};
\node[above,gray] at (-1.5,-3) {\tiny $(z\bar{*}R_3(x,y))*x$};
\node[left] at (-2.55,-4) {\tiny $\phi(z\bar{*}R_3(x,y),x)$};
\node[left] at (-2.55,-6) {\tiny $-\phi(z\bar{*}R_3(x,y),R_3(x,y))$};
\node[right] at (-.9,-5.1) {\tiny $\phi_2(x,y)$};
\node[below,gray] at (4,-7) {\tiny $z$};
\node[above,gray] at (2,-3) {\tiny $x$};
\node[below,gray] at (6,-7) {\tiny $y$};
\node[right,gray] at (5.7,-5) {\tiny $z\bar{*}y$};
\node[above,gray] at (4,-3) {\tiny $(z\bar{*}y)*R_3(y,x)$};
\node[below,gray] at (2,-7) {\tiny $R_3(x,y)$};
\node[above,gray] at (6,-3) {\tiny $R_3(y,x)$};
\node[right] at (5.3,-4) {\tiny $\phi(z\bar{*}y,R_3(y,x))$};
\node[right] at (5.5,-6) {\tiny $-\phi(z\bar{*}y,y)$};
\node[left] at (3.3,-5.1) {\tiny $\phi_2(x,y)$};
\end{tikzpicture}
\vspace{.2in}
		\caption{Bondle coloring and Boltzmann weight for $bRIIIb$.}
		\label{weightax2}
\end{figure}

From Figure~\ref{weightax2} we obtain the following condition:
\begin{equation}\label{bondle2}
\begin{split}
\phi(z \bar{*} R_3(x,y),x)- \phi(z \bar{*} R_3(x,y),R_3(x,y))+ \phi_2(x,y) = \\\phi(z \bar{*} y, R_3(y,x)) - \phi(z \bar{*}y,y) + \phi_2(x,y)\quad \text{for all $x,y,z \in X$}.  
\end{split}
\end{equation}

\begin{figure}[h]
\begin{tikzpicture}[use Hobby shortcut,add arrow/.style={postaction={decorate}, decoration={
  markings,
  mark=at position 1 with {\arrow[scale=1,>=stealth]{>}}}}]
\begin{knot}
 \strand (-2,-11).. (-1.2,-10)..(-1.5,-9)..(-1,-8.5)..(0,-8)[add arrow];
 \strand (-2,-8)..(-1,-8.5)..(-.5,-9)..(-.8,-10)..(0.3,-11)[add arrow];

 \draw (-1.2,-10.3) to (-.8,-10.3);
 \draw (-1.25,-10.5) to (-.75,-10.5);
 
 \draw[<->] (.5,-9.5)..(1.5,-9.5);
 
 \strand (1.9,-11)..(3,-10.5)..(3.5,-10)..(3.2,-9)..(4.1,-8)[add arrow];
 \strand (2,-8)..(2.8,-9)..(2.5,-10)..(3,-10.5)..(4,-11)[add arrow];

 \draw (2.8,-8.9) to (3.2,-8.9);
 \draw (2.8,-8.7) to (3.2,-8.7);
 
\end{knot}
\node[above,gray] at (-2,-8) {\tiny $x$};
\node[below,gray] at (-2,-11) {\tiny $y$};
\node[above,gray] at (0,-8) {\tiny $R_3(y,x)$};
\node[right,gray] at (-.5,-10) {\tiny $x \bar{*} R_3(y,x)$};
\node[left] at (-1.5,-10.4) {\tiny $\phi_2(y, x \bar{*}R_3(y,x))$};
\node[left] at (-1.5,-8.5) {\tiny $-\phi(x \bar{*} R_3(y,x), R_3(y,x))$};

\node[above,gray] at (2,-8) {\tiny $x$};
\node[below,gray] at (2,-11) {\tiny $y$};
\node[above,gray] at (4,-8) {\tiny $R_3(y,x)$};
\node[left,gray] at (2.6,-9.1) {\tiny $R_3(x,y)$};
\node[right] at (3.5,-10.4) {\tiny $-\phi(R_3(x,y)\bar{*}y,y)$};
\node[right] at (3.5,-9) {\tiny $\phi_2(x,y)$};
\end{tikzpicture}
\vspace{.2in}
		\caption{Bondle coloring and Boltzmann weight for $bRII$.}
		\label{weightax3}
\end{figure}
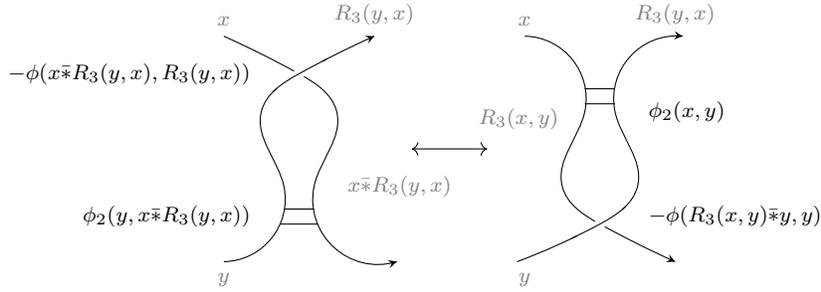

Lastly, Figure~\ref{weightax3} imposes the following condition:
\begin{equation}\label{bondle3}
\begin{split}
-\phi(x \bar{*} R_3(y,x), R_3(y,x))+\phi_2(y, x \bar{*}R_3(y,x)) = \\ \phi_2(x,y) - \phi(R_3(x,y)\bar{*}y,y)\quad \text{for all $x,y,z \in X$}. 
\end{split}
\end{equation}

We now define the extension of the Boltzmann weight invariant to folded molecular chain models. 

\begin{definition}\label{statesum}
Let $(X, *, R_1,R_2,R_3)$ be an oriented bondle and let $A$ be an abelian group. Assume that $\phi,\phi_1,\phi_2: X \times X \rightarrow A$ satisfy the equations (\ref{cocycle1})-(\ref{bondle3}). Then the \emph{state sum} for a folded molecular chain diagram $D$ of a folded molecular chain model $P$ is given by
\[\Phi_X^{\phi,\phi_1,\phi_2}(D) =\sum_{\mathcal{C}}\prod_{\tau} \psi(x,y),\]
where the product is taken over all classical crossings and bonds, denoted by $\tau$, of the diagram $D$. Furthermore, the sum is taken over all possible colorings, denoted by $\mathcal{C}$, of the diagram $D$.
\end{definition}

Note that in the above definition $\psi (x,y) = \phi(x,y)^{\pm 1}$ at classical positive and negative crossings. Additionally, $\psi(x,y) = \phi_1(x,y)$ at a bond with parallel strands and $\psi(x,y) = \phi_2(x,y)$ at a bond with anti-parallel strands. Before we compute explicit examples of the Boltzmann weight enhanced invariant, we first prove that the state sum, $\Phi_X^{\phi,\phi_1,\phi_2}$, is an invariant of folded molecular chain models.

\begin{theorem}
Let $\phi, \phi_1,\phi_2: X \times X \rightarrow A$ be maps satisfying the conditions in Definition~\ref{statesum}.  The state sum associated with the functions $\phi, \phi_1$ and $\phi_2$ is invariant under the moves in Figures~\ref{weightax1}, \ref{weightax2} and \ref{weightax3}, so it defines an invariant of folded molecular chains with bonds. 
\end{theorem}
\begin{proof}
The proof is similar to that in classical knot theory.   There is one-to-one correspondence between colorings before and after each of the generating set of Reidemeister moves involving bonds.  The equations~(\ref{sing1}), (\ref{sing2}), (\ref{sing3}), (\ref{bondle1}), (\ref{bondle2}) and (\ref{bondle3}) imply the result.  
\end{proof}
The following proposition gives explicit formulas of the Boltzmann weights coming from $\phi, \phi_1$ and $\phi_2$.
\begin{proposition}\label{boltzmann}
Let $(X, *, R_1, R_2,R_3)$ be an oriented bondle and $\phi, \phi_1, \phi_2 : X \times X \rightarrow \mathbb{Z}_n$. If $\phi(x,y)=0$, $\phi_1(x,y)=a$, and $\phi_2(x,y)=b$ for all $x,y  \in X$ and $a,b \in \mathbb{Z}_n$. Then $\phi$, $\phi_1$, and $\phi_2$ are Boltzmann weights. 
\end{proposition}

\begin{proof}
We verify the Boltzmann weight axioms (\ref{cocycle1})-(\ref{bondle3}). Note that equations (\ref{cocycle1}), (\ref{cocycle2}), (\ref{sing2}), and  ($\ref{bondle2}$) are automatically satisfied since $\phi(x,y)=0$ for all $x,y \in X$. The axioms (\ref{sing1}) and (\ref{sing3}) reduce to $a=a$. Lastly, axioms (\ref{bondle1}) and (\ref{bondle3}) reduce to $b = b$. Therefore, $\phi$, $\phi_1$, and $\phi_2$ are Boltzmann weight of $X$ on $\mathbb{Z}_n$.
\end{proof}

Noticing that the target group of the maps $\phi, \phi_1$ and $\phi_2$ is a finite cyclic group, $\mathbb{Z}_n = \langle u, u^n = 1 \rangle$, then the state sum invariant will have the form $\Phi_X^{\phi, \phi_1, \phi_2}(D) = \sum_{i=0}^{n-1}a_i u^i$.

\begin{figure}[h]
    \centering
    \includegraphics[scale=2.2]{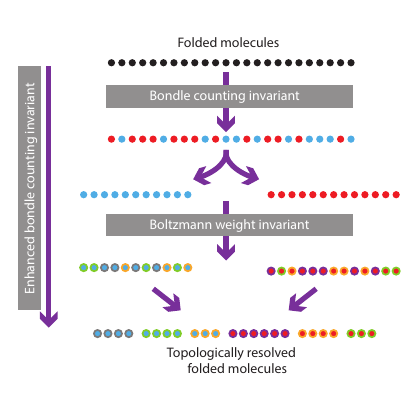}
    \caption{Schematic of the coloring enhancement.}
    \label{schematic}
\end{figure}

Before applying the state sum invariant using Boltzmann weights to examples, we present an overview of the 
invariant's implementation.  We also provide a visual representation of the bondle counting invariant implementation and the state sum invariant using Boltzmann weights, 
see Figure~\ref{schematic}. Suppose one wants to do a systematic study of distinguishing folded molecular chains. In that case, one starts from a collection (or database) of folded molecular chains, illustrated. The first thing to do is to apply the bondle counting invariant to make clusters of folded molecular chains with the same number of colorings.  Now within a fixed grouping, the bondle counting invariant becomes useless in distinguishing pairs of folded molecules.  To get a more refined clustering within each cluster, we can compute the state sum invariant of each folded molecular chain using the Boltzmann weights.

\section{Examples}\label{Examples}
In this section, we compute the state sum invariant using Boltzmann weights of projections of example proteins. We show that the 
new invariant is a proper enhancement of the bondle coloring invariant by presenting diagrams of proteins with the same bondle coloring invariant but with different Boltzmann weights. Therefore, the state sum invariant using Boltzmann weights provides a finer invariant than the coloring invariant for the classification of proteins into distinct topological types. The following bondles and Boltzmann weights were obtained by Proposition~\ref{bondleEx} and Proposition~\ref{boltzmann}. In addition, all of the computations were performed using \emph{Mathematica software} and confirmed independently by \emph{Maple software}.

\begin{example}
Let $(X,*, R_1,R_2,R_3)$ be an oriented bondle with operations defined on $X=\Z_{15}$ by, $x* y = 4x-3y = x \bar{*}y$, $R_1(x,y) = 3x-2y$, $R_2(x,y) = -8x+9y$, and $R_3(x,y) = 6x+5y$. The two proteins $P_1$ and $P_2$ pictured below both have $45$ colorings by the oriented bondle $X$.
\begin{figure}[h]
\begin{tikzpicture}[use Hobby shortcut,scale=1.4]
\begin{knot}[
consider self intersections=true,
  ignore endpoint intersections=false,
  clip width=4,
  flip crossing/.list={2}
]
\strand[decoration={markings,mark=at position 6cm with
    {\arrow[scale=1.5,>=stealth]{>}}},postaction={decorate},decoration={markings,mark=at position 24cm with{\arrow[scale=1.5,>=stealth]{>}}},postaction={decorate}] (-1,1)..(-.5,1)..(0,1) ..(0,-2)..(2,-3)..(4,-2)..(3,-1)..(4,0)..(3,.5)..(2.5,.5)..(1,-1)..(2,-2)..(4,-1)..(3,0)..(4,1.5)..(2,1.5)..(0,2)..(-1,2); 

    \draw (2,-3) to (2,-2);
    \draw (2.1,-3) to (2.1,-2.05);
    \draw (-.7,1) to (-.7,2.05);
    \draw (-.6,1) to (-.6,2.05);
    \draw (3,-1.05) to (4.02,-1.05);
    \draw (3,-1.15) to (4.05,-1.15);
    

\end{knot}
\node[below] at (-1,1) {\tiny $v$};
\node[left] at (-.1,-1) {\tiny $x$};
\node[below] at (3,-3.2) {\tiny $R_1(x,y)$};
\node[above] at (3,-2.1) {\tiny $R_2(x,y)$};
\node[right] at (4,-1.5) {\tiny $z$};
\node[right] at (4,-.8) {\tiny $R_1(z,R_1(x,y))$};
\node[left] at (3,-.8) {\tiny $R_2(z,R_1(x,y))$};
\node[right] at (4,0) {\tiny $y$};
\node[above] at (2,1.6) {\tiny $u$};
\end{tikzpicture}
\caption{Diagram of protein $P_1$}
\end{figure}
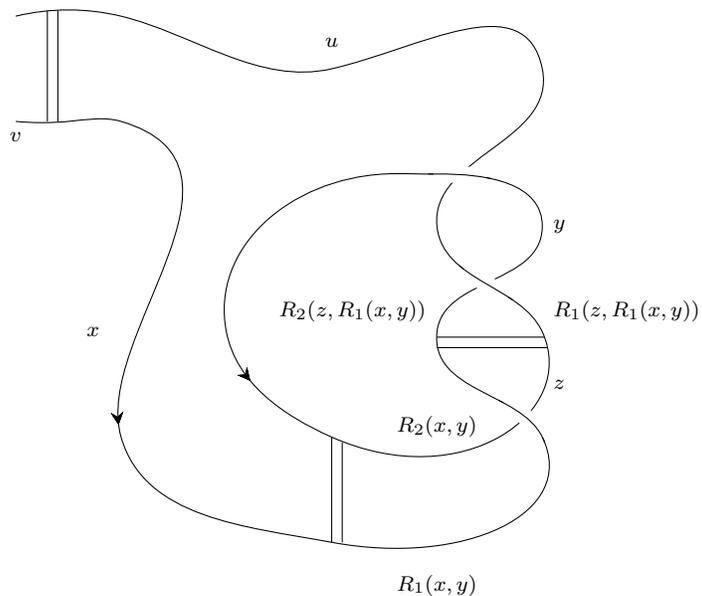

\begin{figure}[h]
\begin{tikzpicture}[use Hobby shortcut,scale=1.4]
\begin{knot}[
consider self intersections=true,
clip width=4,
  ignore endpoint intersections=false,
flip crossing/.list={1,3}
]
\strand[decoration={markings,mark=at position .1 with
    {\arrow[scale=1.5,>=stealth]{>}}},postaction={decorate}] (-1,1)..(-.5,1)..(0,1) ..(0,-2)..(2,-3)..(4,-2)..(3,-1)..(4,0)..(3,.5)..(1.5,.5)..(1,0)..(2.5,-1.5)..(1.5,-2)..(1.3,-1.3)..(4,-1)..(3,0)..(4,1.5)..(2,1.5)..(0,2)..(-1,2); 

    \draw (1.7,-2.95) to (1.7,-2.1);
    \draw (1.8,-2.96) to (1.8,-2.1);
    \draw (-.7,1) to (-.7,2.05);
    \draw (-.6,1) to (-.6,2.05);
    \draw (3,-1.05) to (4.01,-1.05);
    \draw (3,-1.15) to (4.02,-1.15);

\end{knot}
\node[below] at (-1,1) {\tiny $v$};
\node[left] at (-.1,-1) {\tiny $x$};
\node[below] at (3,-3.2) {\tiny $R_3(x,y)$};
\node[right] at (4,-1.5) {\tiny $z$};
\node[right] at (4,-.8) {\tiny $R_1(z,R_3(x,y)*z)$};
\node[right] at (4,0) {\tiny $R_2(z,R_3(x,y)*z)$};
\node[left] at (1.2,-1.5) {\tiny $R_3(y,x)$};
\node[right] at (1.8,.8) {\tiny $y$};
\node[above] at (2,1.6) {\tiny $u$};
\end{tikzpicture}
\caption{Diagram of protein $P_2$}
\end{figure}
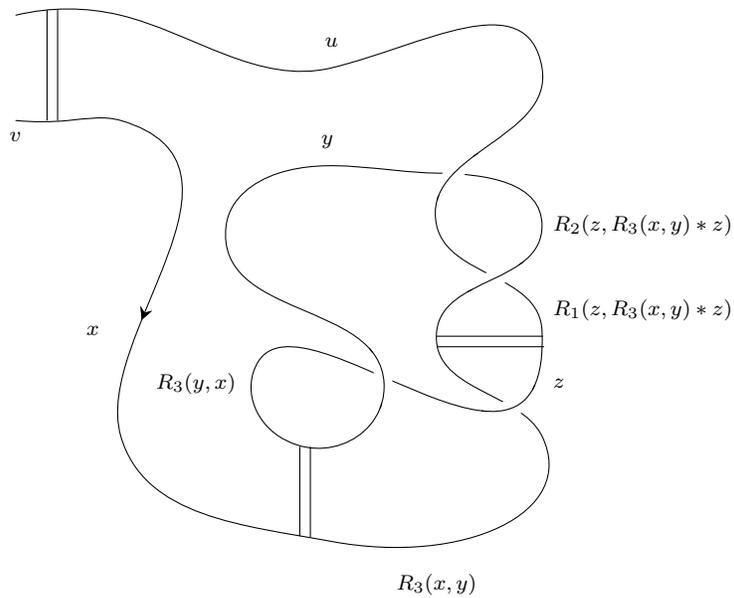

Consider the following Boltzmann weights, 
$\phi,\phi_1, \phi_2: X \times X \rightarrow \Z_6$ with outputs in the set $\Z_6$ of integers modulo $6$ and defined by $\phi(x,y) = 0$, $\phi_1(x,y) = 4$, and $\phi_2(x,y) = 5$. Protein $P_1$ has Boltzmann weight 
\[ \Phi_X^{\phi,\phi_1,\phi_2}(P_1) = 45u. \]
On the other hand, protein $P_2$ has Boltzmann weight
\[  \Phi_X^{\phi,\phi_1,\phi_2}(P_2) = 45u^3.\]
Therefore, $P_1$ and $P_2$ have the same number of colorings, but are distinguished by their Boltzmann weight. 
\end{example}

\begin{example}
Let $(X,*, R_1,R_2,R_3)$ be an oriented bondle with operations defined on $X=\Z_{15}$ by $x* y = 2x-y$, $x\bar{*}y = 8x-7y$, $R_1(x,y) = 4x-3y$, $R_2(x,y) = -6x+7y$, and $R_3(x,y) = 10x-4y$. The two proteins $P_3$ and $P_4$ pictured below both have $75$ colorings by the oriented bondle $X$.

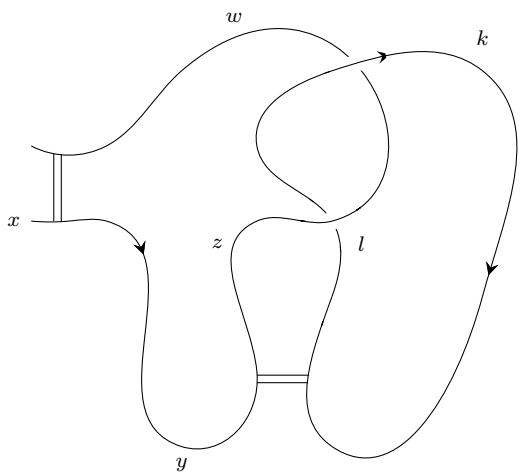
\begin{figure}[h]
\begin{tikzpicture}[use Hobby shortcut,scale=1]
\begin{knot}[
consider self intersections=true,
  ignore endpoint intersections=false,
  clip width=4,
flip crossing/.list={1,3}
]
\strand[decoration={markings,mark=at position .1 with
    {\arrow[scale=1.5,>=stealth]{>}}},postaction={decorate}] (-1,1)..(-.5,1)..(0,1)..(1,-2)..(2,-1)..(2,1)..(3,1)..(1,3)..(0,2)..(-1,2); 
\strand[decoration={markings,mark=at position .4 with
    {\arrow[scale=1.5,>=stealth]{<}}},postaction={decorate}] ([closed]3,0)..(3,-2)..(5,0)..(5,3)..(3,3)..(2,2)..(3,1);
    \draw (-.7,1) to (-.7,1.9);
    \draw (-.6,1) to (-.6,1.89);
    \draw (2,-1.05) to (2.69,-1.05);
    \draw (2,-1.15) to (2.68,-1.15);
\end{knot}
\node[left] at (-1,1) {\tiny $x$};
\node[below] at (1,-2) {\tiny $y$};
\node[left] at (1.7,.7) {\tiny $z$};
\node[above] at (1.7,3.5) {\tiny $w$};
\node[right] at (3.2,.7) {\tiny $l$};
\node[above] at (5,3.2) {\tiny $k$};
\end{tikzpicture}
\caption{Diagram of protein $P_3$}
\end{figure}

\begin{figure}[h]
\begin{tikzpicture}[use Hobby shortcut,scale=1]
\begin{knot}[
consider self intersections=true,
clip width=4,
  ignore endpoint intersections=false,
flip crossing/.list={1,3}
]
\strand[decoration={markings,mark=at position .1 with
    {\arrow[scale=1.5,>=stealth]{>}}},postaction={decorate}] (-1,1)..(-.5,1)..(0,1)..(1,-2)..(2,-1)..(2,1)..(3,1)..(1,3)..(0,2)..(-1,2); 
\strand[decoration={markings,mark=at position .4 with
    {\arrow[scale=1.5,>=stealth]{>}}},postaction={decorate}] ([closed]3,0)..(3,-2)..(5,0)..(5,3)..(3,3)..(2,2)..(3,1);
    \draw (-.7,1) to (-.7,1.9);
    \draw (-.6,1) to (-.6,1.89);
    \draw (2,-1.05) to (2.69,-1.05);
    \draw (2,-1.15) to (2.68,-1.15);
\end{knot}

\node[left] at (-1,1) {\tiny $x$};
\node[below] at (1,-2) {\tiny $y$};
\node[left] at (1.7,.7) {\tiny $z$};
\node[above] at (1.7,3.5) {\tiny $w$};
\node[right] at (3.2,.7) {\tiny $l$};
\node[above] at (5,3.2) {\tiny $k$};
\end{tikzpicture}
\caption{Diagram of protein $P_4$}
\end{figure}
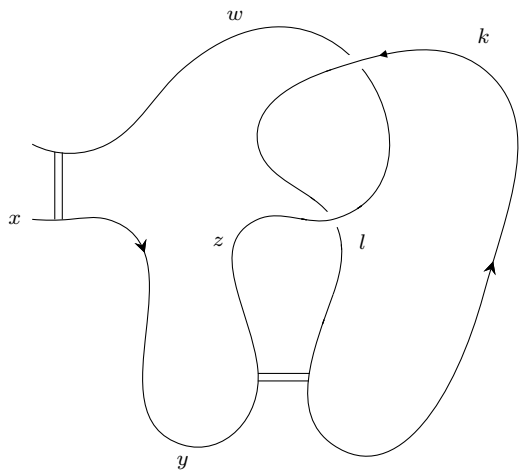

Consider the following Boltzmann weights, 
$\phi,\phi_1, \phi_2: X \times X \rightarrow \Z_6$ with outputs in the set $\Z_6$ of integers modulo $6$ and defined by $\phi(x,y) = 0$, $\phi_1(x,y) = 1$, and $\phi_2(x,y) = 3$. Protein $P_3$ has Boltzmann weight 
\[ \Phi_X^{\phi,\phi_1,\phi_2}(P_3) = 75. \]
On the other hand, protein $P_4$ has Boltzmann weight
\[  \Phi_X^{\phi,\phi_1,\phi_2}(P_4) = 75u^2.\]
Therefore, $P_3$ and $P_4$ have the same number of colorings, but are distinguished by their Boltzmann weight. 
\end{example}
\section{Conclusion}
We developed an enhanced coloring scheme that enables the topological classification of folded linear chains. The approach applies to any given folded polymer chains, featuring arbitrary arrangements of intra-chain contacts and chain crossings. The enhancement of the bondle coloring invariant is especially interesting since, for each oriented bondle, there are several Boltzmann weight choices. By selecting an appropriate Boltzmann weight, we can emphasize the intra-chain contacts with parallel strands, intra-chain contacts with anti-parallel strands, or the chain crossings.  In \cite{ADEM}, the Gauss code of folded linear chains diagrams was introduced; this allows us to input a folded linear chain diagram into a computer for computations. Therefore, with the use of Gauss codes, we can compute the bondle coloring invariant and the Boltzmann weight enhancement for a large number of folded linear chain diagrams. The only limitation would be computing power, but we can also select a smaller bondles set and a smaller abelian group to decrease the computation complexity. Since our enhancement is built on the idea of an oriented bondle, we are restricted to segmentations of $\beta$-pleated sheets. This means that the following three diagrams are still indistinguishable. 

\begin{figure}[h]
    \centering
    \includegraphics[scale=.6]{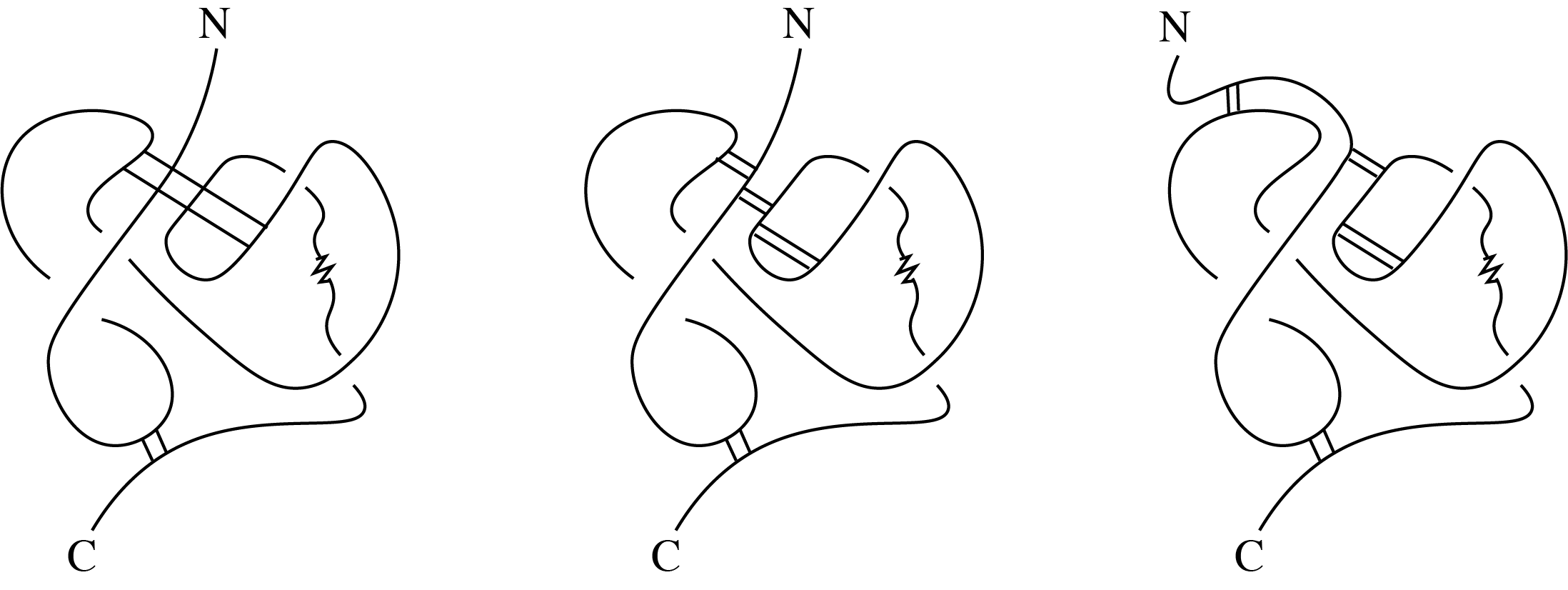}
    \caption{Three indistinguishable protein models}
    \label{fig:my_label}
\end{figure}
\medskip

\section*{Acknowledgement} The authors would like to thank Colin Adams for fruitful discussions which improved the paper.\\
M.E. was partially supported by Simons Foundation collaboration grant 712462.

\medskip


\end{document}